\documentclass[a4paper,11pt,leqno,twoside]{article}
\usepackage{amsmath,amsthm,amssymb}
\usepackage[all]{xy}
\numberwithin{equation}{section}
\usepackage[dvips]{graphicx,psfrag}
\usepackage[usenames]{color}

\input{mymacros}
\newcommand{\limup}{\lim_{q \uparrow 1}}
\newcommand{\limdo}{\lim_{q \downarrow 0}}
\newcommand{\sumn}{\sum_{n \in \bZ\backslash\{0\}}}
\newcommand{\p}{\partial}
\newcommand{\nn}{^{(\n)}}

\theoremstyle{plain}

\newtheorem{thm}{Theorem}[section]
\newtheorem{prop}[thm]{Proposition}
\newtheorem{lem}[thm]{Lemma}
\newtheorem{cor}[thm]{Corollary}

\theoremstyle{definition}

\newtheorem{remark}[thm]{Remark}

\title{$q$-Analogues of the Riemann zeta, the Dirichlet $L$-functions,\\ 
and a crystal zeta function} 
\author{Kenichi Kawagoe, Masato Wakayama and Yoshinori Yamasaki}

\begin{document}

\setlength{\baselineskip}{16pt}
\maketitle

\begin{abstract}
 A $q$-analogue $\z_q(s)$ of the Riemann zeta function $\z(s)$ was
 studied in \cite{KanekoKurokawaWakayama2003} via a certain $q$-series
 of two variables. We introduce in a similar way a $q$-analogue of the
 Dirichlet $L$-functions and make a detailed study  of them, including
 some issues concerning the classical limit of $\z_q(s)$ left open in
 \cite{KanekoKurokawaWakayama2003}. We also examine a ``crystal'' limit 
 (i.e. $q \downarrow 0$) behavior of $\z_q(s)$. The $q$-trajectories of
 the trivial and essential zeros of $\z(s)$ are investigated numerically
 when $q$ moves in $(0,1]$. Moreover, conjectures for the crystal limit
 behavior of zeros of $\z_q(s)$ are given.
\end{abstract}

%
\section{Introduction}
%

 There are fairly plenty of possibilities for defining $q$-analogues of
 the Riemann zeta function $\z(s)=\sum^{\infty}_{n=1}n^{-s}$ in the
 convergent region $\Re(s)>1$, (See, e.g., \cite{Satoh1989},
 \cite{Cherednik2001} and \cite{KanekoKurokawaWakayama2003}). Among
 them, in \cite{KanekoKurokawaWakayama2003} Kaneko, Kurokawa and the
 second author introduced a certain $q$-analogue $\z_q(s)$ of $\z(s)$
 which is meromorphically extended to the entire plane $\bC$ and indeed
 gives a proper $q$-analogue in the sense that the classical limit of
 $\z_q(s)$ exists and equals $\z(s)$ for {\it all} $s \in \bC$. We
 briefly recall the story. Suppose that $0<q<1$. Let $f_q(s,t)$ be a
 function of the two complex variables $s$ and $t$ defined by the series
 $f_q(s,t):=\sum^{\infty}_{n=1}q^{nt}[n]^{-s}_q$, where
 $[n]_q:=(1-q^n)/(1-q)$ denotes a $q$-analogue of the number $n$. It is
 clear that the series converges absolutely for $\Re(t)>0$. Obviously,
 $\lim_{q \uparrow 1}f_q(s,t)=\z(s)$ holds in the convergent region
 $\Re(t)>0$. Among these $f_q(s,t)$, the $q$-analogue $\z_q(s)$ of
 $\z(s)$ is defined by setting $\z_q(s):=f_q(s,s-1)$ for
 $\Re(s)>1$. Then it was shown that $\z_q(s)$ can be meromorphically
 continued to $\bC$ and moreover that $\limup\z_q(s)=\z(s)$ for all
 $s\in \bC$. Note that, however, neither an Euler product expression nor
 a functional equation can be expected for $\z_q(s)$.  

 The initial aim of the present paper is to generalize the result on
 $\z_q(s)$ to the cases of the Dirichlet $L$-functions and to make a
 much detailed study of those $q$-analogues including some results about
 the case of the Riemann zeta function not even clarified in
 \cite{KanekoKurokawaWakayama2003}. Let $N \in \bN$ and $\c$ be a
 Dirichlet character modulo $N$. Let
 $L(s,\c)=\sum^{\infty}_{n=1}\c(n)n^{-s}$ $(\Re(s)>1)$ be the Dirichlet
 $L$-functions. To explain the result precisely, let us define a series
 $f_q(s,t,\c)$ similar to $f_q(s,t)$ by the formula
 $f_q(s,t,\c)=\sum^{\infty}_{n=1}\c(n)q^{nt}[n]^{-s}_q$ for
 $\Re(t)>0$. It is then proved that $f_q(s,t,\c)$ can be meromorphically
 continued to the entire plane $\bC$. Putting $L_q(s,\c)=f_q(s,s-1,\c)$,
 as one may expect from the study of $\z_q(s)$, we actually prove that
 $\limup L_q(s,\c)=L(s,\c)$ for all $s\in\bC$. Moreover, not only do we
 treat the case of $L_q(s,\c)$, but also we show that each of the
 functions $f_q(s,s-\n,\c)$, $\n=2,3,\ldots$, gives a proper
 $q$-analogue of $L(s,\c)$, and, what is more, that only these functions
 $f_q(s,s-\n,\c)$ $(\n\in\bN)$ can realize such true $q$-analogues of
 $L(s,\c)$ in the family of the meromorphic functions of the form
 $f_q(s,\varphi(s),\c)$ provided $\varphi(s)$ is a non-constant
 meromorphic function on $\bC$. If $\c$ is not principal, however, in
 addition to the case $\varphi(s)=s-\n$ for $\n\in\bN$, the constant
 function $\varphi(s)=\m\in\bN$ gives also a true $q$-analogue of
 $L(s,\c)$. Our analysis is based on the use of a $q$-analogue of the
 Hurwitz zeta function defined similarly.     

 A numerical analysis of the zeros of $\z_q(s)$, which is the second
 purpose of this paper, is developed in the last section $\S 3$. We
 first examine a ``crystal'' limit (i.e. the pointwise limit for
 $q \downarrow 0$) of $\z_q(s)$ and show that an analogue of the Riemann
 hypothesis holds for such a crystal Riemann zeta function. Furthermore,
 the $q$-trajectories of the trivial zeros and the essential zeros of
 $\z(s)$ are numerically investigated, that is, the zeros of $\z_q(s)$
 are studied with Maple 8 \cite{maple}. Especially, we observe that the
 limit point of the $q$-trajectory of each essential zero of $\z(s)$
 falls on either $0$ or the one of the trajectory of the some trivial
 zero of $\z(s)$, i.e. on a negative integer point. We then give some
 conjectures concerning the crystal limit behavior of zeros of
 $\z_q(s)$.

%
\section{On classical limits}
%

We study the $q$-analogues of the Dirichlet $L$-function
$L(s,\c)=\sum^{\infty}_{n=1}\c(n)n^{-s}$ and the Hurwitz zeta function
$\z(s,a)=\sum^{\infty}_{n=0}(n+a)^{-s}$ defined respectively by the
series    
\[
 f_q(s,t,\c):=\sum^{\infty}_{n=1}\frac{\c(n)q^{nt}}{[n]_q^s}
    \quad \textrm{and} \quad 
  g_q(s,t,a):=\sum^{\infty}_{n=0}\frac{q^{(n+a)t}}{[n+a]_q^s},
\]
where $\c$ is a Dirichlet character modulo $N$ and $0<a\le1$. These
series converge absolutely for $\Re(t)>0$. Obviously,
$f_q(s,t)=f_q(s,t,{\bf 1})=g_q(s,t,1)$ ({\bf 1} denotes the trivial
character), where $f_q(s,t)$ is the function discussed in
\cite{KanekoKurokawaWakayama2003} (See \S 1). We put
$L\nn_q(s,\c):=f_q(s,s-\n,\c)$, $\z\nn_q(s):=f_q(s,s-\n,{\bf 1})$ for
$\n \in \bN$, and $\z_q(s)=\z^{(1)}_q(s)$. These series converge
absolutely for $\Re(s)>\n$. We will show later that if $\c$ is not the 
principal character, $L\nn_q(s,\c)$ is holomorphic at
$s=\n,\n-1,\n-2,\ldots,1$, and if $\c$ is the principal character,
$L\nn_q(s,\c)$ has simple poles at these points. We find, in particular,
that $\z\nn_q(s)$ has simple poles at these points (See
Proposition\,\ref{prop:specialvalue}). Further we define an entire
function $L^{\m}_q(s,\c)$ for $\m\in\bN$ by 
\[
 L^{\m}_q(s,\c)
:=f_q(s,\m,\c)=\sum^{\infty}_{n=1}\frac{\c(n)q^{{\m}n}}{[n]_q^s}.
\]   
 
 We first prove the following theorem, which, in particular, gives an
 affirmative answer of the question in
 \cite{KanekoKurokawaWakayama2003}, Remark $(1)$ on p.\,$179$. 


\begin{thm}\label{thm:hz}
 Let $t=\varphi(s)$ be a meromorphic function on $\bC$. Then the formula 
\[
 \limup g_q(s,\varphi(s),a)=\z(s,a) \quad (s \in \bC)
\]
holds if and only if the function $\varphi(s)$ can be written as
 $\varphi(s)=s-\n$ for some $\n \in \bN$. Namely, $g_q(s,s-\n,a)$ gives
 a true $q$-analogue of $\z(s,a)$ and these are the only true
 $q$-analogues among the functions of the form $g_q(s,\varphi(s),\c)$.   
\end{thm}

To prove Theorem\,\ref{thm:hz}, recall the Euler-Maclaurin summation
formula (see, \cite{Titchmarsh1986}). For integers $b, c$ satisfying
$b<c$, a $C^{\infty}$-function $f(x)$ on $[b,\infty)$, and an arbitrary
integer $M \ge 0$, we have    
\begin{align}\label{for:EM} 
  \sum^c_{n=b}f(n)
&=\int^c_b
 f(x)dx+\frac{1}{2}(f(b)+f(c))+\sum^M_{l=1}\frac{B_{l+1}}{(l+1)!}(f^{(l)}(c)-f^{(l)}(b))\\    
&\ \ \ -\frac{(-1)^{M+1}}{(M+1)!}\int^c_b\widetilde{B}_{M+1}(x)f^{(M+1)}(x)dx,\nonumber
\end{align}
where $B_{n}$ is the Bernoulli number and $\widetilde{B}_{M+1}(x)$ is
the periodic Bernoulli polynomial defined by
$\widetilde{B}_{k}(x)=B_{k}(x-[x])$ with $[x]$ being the largest integer
not exceeding $x$. Assume $\Re(t)>0$. Since
$g_q(s,t,a)=(1-q)^s\sum^{\infty}_{n=0}q^{(n+a)t}(1-q^{n+a})^{-s}$, the
Euler-Maclaurin summation formula \eqref{for:EM} gives  
\begin{align*}
  g_q(s,t,a)
&=\frac{1}{2}q^{at}\Bigl(\frac{1-q^a}{1-q}\Bigr)^{-s}
-\frac{1}{12}q^{at}\Bigl(\frac{1-q^a}{1-q}\Bigr)^{-s}\Biggl(\Bigl(\frac{1-q^a}{\log{q}}\Bigr)^{-1}s+(t-s)\log{q}\Biggr)\\
&\ \ \ +(1-q)^sI^0_0(s,t;q)+(1-q)^s\sum^2_{\e=0}a_{\e}(s,t;q)I_{\e}(s,t;q),  
\end{align*}
where
\begin{align*}
  I^0_0(s,t;q)&=\int^{\infty}_0q^{(x+a)t}(1-q^{x+a})^{-s}dx,\\
 I_{\e}(s,t;q)&=\int^{\infty}_0\widetilde{B}_2(x)q^{(x+a)t}(1-q^{x+a})^{-s-\e}dx 
 \quad (\e=0,1,2),\\  
 a_{\e}(s,t;q)&=
    \begin{cases}
     -\frac{1}{2}(\log{q})^2(s-t)^2 \quad & (\e=0),\\
      \frac{1}{2}(\log{q})^2s(2s-2t+1) \quad & (\e=1),\\
     -\frac{1}{2}(\log{q})^2s(s+1) \quad & (\e=2).
    \end{cases}
\end{align*}
Change the variables $u=q^{x+a}$ in the first integral
$I^0_0(s,t;q)$. Then    
\begin{align}\label{for:00} 
 I^0_0(s,t;q)=-\frac{1}{\log{q}}b_{q^a}(t,-s+1) \qquad
 \Bigl(\ \textrm{or}\quad b_{q^a}(\a,\b)=-\log{q}\cdot
 I^0_0(-\b+1,\a;q)\Bigr),
\end{align}
where $b_{q}(\a,\b)$ is the incomplete beta function defined by 
\[
 b_{q}(\a,\b)=\int^q_0u^{\a-1}(1-u)^{\b-1}du \qquad (0<q<1).
\]
The integral $b_{q}(\a,\b)$ converges absolutely for
$\Re(\a)>0$. Moreover, the Fourier expansion  
\begin{equation}\label{for:Fourier}
 \widetilde{B}_k(x)=-k!\sumn\frac{e^{2\pi inx}}{(2\pi in)^k}
\end{equation}
gives 
\begin{equation}\label{for:I123}
 I_{\e}(s,t;q)=\frac{2}{\log{q}}\sumn\frac{e^{-2\pi ina}}{(2\pi
 in)^2}b_{q^a}(t+\d n ,-s-\e+1),
\end{equation}
where $\d:=2\pi i/\log{q}$. Note that since $b_{q^a}(t+\d n ,-s-\e+1)$ 
converges absolutely for $\Re(t)>0$ and is uniformly bounded with
respect to $n$, the sum \eqref{for:I123} converges absolutely. Take $M 
\ge 2$ satisfying $\Re(s)>-M$. Then, integral by parts yields
\begin{align}\label{for:ac-beta}
  b_q(\a,\b)
&=\sum^{M-1}_{l=1}(-1)^{l-1}\frac{(1-\b)_{l-1}}{(\a)_l}q^{\a+l-1}(1-q)^{\b-l}\\ 
&\ \ \ +(-1)^{M-1}\frac{(1-\b)_{M-1}}{(\a)_{M-1}}b_q(\a+M-1,\b-M+1),\nonumber
\end{align}
where $(s)_k:=s(s+1)(s+2)\cdots(s+k-1)$. Hence, by \eqref{for:ac-beta},
we have    
\begin{multline}\label{for:beta1} 
 \ \ b_{q^a}(t+\d n ,-s-\e+1)
= e^{2\pi
 ina}\Biggl\{(1-q^a)^{-s-\e+1}\sum^{M-1}_{k=1}\frac{(-1)^{k-1}(s+\e)_{k-1}}{(1-q^a)^k(t+\d n)_k}q^{a(t+k-1)}\Biggr.\\  
\ \ \ \ \ \ \ \  \ \ \ \ \ \ \  \Biggl.-\log{q}(1-q)^{-s-\e}\frac{(-1)^{M-1}(s+\e)_{M-1}}{(1-q)^{M-1}(t+\d 
 n)_{M-1}}\int^{\infty}_0 e^{2\pi inx}q^{(x+a)(t+M-1)}\Bigl(\frac{1-q^{x+a}}{1-q}\Bigr)^{-s-\e+1-M}dx\Biggr\}.
\end{multline}
This gives an analytic continuation of $b_{q^a}(t+\d n ,-s-\e+1)$ to the
region $\Re(t)>1-M$. Hence, by \eqref{for:beta1}, we obtain   
\begin{equation}\label{for:C123}
 (1-q)^sa_{\e}(s,t;q)I_{\e}(s,t;q)=\frac{2(1-q^a)^{-\e+1}}{\log{q}}a_{\e}(s,t;q)C_{\e}(s,t,M;q)=:\widetilde{C}_{\e}(s,t,M;q).
\end{equation}
Here we put
\begin{multline}\label{for:C123'}
 \ \ C_{\e}(s,t,M;q)=\sumn\frac{1}{(2\pi
 in)^2}\Biggl\{\Bigl(\frac{1-q^a}{1-q}\Bigr)^{-s}\sum^{M-1}_{k=1}\frac{(-1)^{k-1}(s+\e)_{k-1}}{(1-q^a)^k(t+\d n)_k}q^{a(t+k-1)}\Biggr.\\
\ \ \ \Biggl.-\Bigl(\frac{1-q^a}{\log{q}}\Bigr)^{-1}\Bigl(\frac{1-q^a}{1-q}\Bigr)^{\e}\frac{(-1)^{M-1}(s+\e)_{M-1}}{(1-q)^{M-1}(t+\d
 n)_{M-1}}\int^{\infty}_0 e^{2\pi inx}q^{(x+a)(t+M-1)}\Bigl(\frac{1-q^{x+a}}{1-q}\Bigr)^{-s-\e+1-M}dx\Biggr\}.
\end{multline}


\begin{lem}\label{lem:C123}
 Let $M\ge 2$ be an integer satisfying $\Re(s)>-M$. Then $\limup
 \widetilde{C}_{\e}(s,t,M;q)$ $(\e=0,1,2)$ exist and are given by  
\begin{align*}
\limup \widetilde{C}_{0}(s,t,M;q)
&=\limup \widetilde{C}_{1}(s,t,M;q)=0,\\
 \limup \widetilde{C}_{2}(s,t,M;q)
&=\sum^{M}_{k=2}\frac{B_{k+1}}{(k+1)!}(s)_{k}a^{-s-k}-\frac{(s)_{M+1}}{(M+1)!}\int^{\infty}_{0}\widetilde{B}_{M+1}(x)(x+a)^{-s-1-M}dx. 
\end{align*} 
\end{lem}
\begin{proof}
 The results follow from the fact $\limup(1-q^a)^k(t+\d n)_k=(-2\pi
 ina)^k$, \eqref{for:I123} and \eqref{for:C123'}. 
\end{proof}

By \eqref{for:00}, \eqref{for:I123} and \eqref{for:C123}, we obtain the
following proposition.   


\begin{prop}\label{prop:q-EM}
 Let $M\ge 2$ be an integer satisfying $\Re(s)>-M$. Then we have  
\begin{align*}
 g_q(s,t,a)
&=\frac{1}{2}q^{at}\Bigl(\frac{1-q^a}{1-q}\Bigr)^{-s}
 -\frac{1}{12}q^{at}\Bigl(\frac{1-q^a}{1-q}\Bigr)^{-s}\Biggl(\Bigl(\frac{1-q^a}{\log{q}}\Bigr)^{-1}s+(t-s)\log{q}\Biggr)\\&-\frac{(1-q)^s}{\log{q}}b_{q^a}(t,-s+1)+\sum^2_{\e=0}\widetilde{C}_{\e}(s,t,M;q).
\end{align*}
This gives the analytic continuation of $g_q(s,t,a)$ to the region
 $\Re(t)>1-M$. \qed 
\end{prop}


\begin{proof}[Proof of Theorem\,\ref{thm:hz}.]
 We first show the sufficiency, that is, for each $\n\in\bN$, $\limup
 g_q(s,s-\n,a)=\z(s,z)$ for all $s\in\bC$. Initially, we assume
 $t=s-1$. Though the case is referred in
 \cite{KanekoKurokawaWakayama2003}, Remark $(1)$ on p.\,$184$, we give
 a proof for completeness. The condition $\Re(t)>0$ implies
 $\Re(s)>1$. In this case, we can evaluate $b_{q^a}(s-1,-s+1)$ by an
 elementary function. In fact, in general, since 
\begin{equation}\label{for:special-beta}
 b_{q}(\a-\n+1,-\a)=-\sum^{\n-1}_{r=0}\frac{(-\n+1)_r}{(-\a)_{r+1}}q^{\a-\n+1}(1-q)^{-\a+r}
 \qquad (\n\in\bN, \Re(\a)>\n-1), 
\end{equation}
we have $b_{q^a}(s-1,-s+1)=\frac{1}{s-1}q^{a(s-1)}(1-q^a)^{-s+1}$. Therefore,
 by Lemma\,\ref{lem:C123}, for $\Re(s)>2-M$ we obtain    
\begin{align}\label{for:limg_q}
 \limup g_q(s,s-1,a)
&=\frac{1}{s-1}a^{1-s}+\frac{1}{2}a^{-s}+
  \sum^{M}_{k=1}\frac{B_{k+1}}{(k+1)!}(s)_{k}a^{-s-k}\\
&-\frac{(s)_{M+1}}{(M+1)!}\int^{\infty}_{0}\widetilde{B}_{M+1}(x)(x+a)^{-s-1-M}dx.\nonumber 
\end{align}
 Comparing the equation \eqref{for:limg_q} with (the analytic
 continuation of) the Hurwitz zeta function  
\begin{align}\label{for:EM-hz} 
 \z(s,a)
&=\frac{1}{s-1}a^{1-s}+\frac{1}{2}a^{-s}+\sum^M_{l=1}\frac{B_{l+1}}{(l+1)!}(s)_la^{-s-l}\\
&-\frac{(s)_{M+1}}{(M+1)!}\int^{\infty}_0\widetilde{B}_{M+1}(x)(x+a)^{-s-M-1}dx
 \qquad (\Re(s)>-M),\nonumber 
\end{align} 
 we have the claim for $\n=1$.

Suppose next $t=s-\n$ for $\n\ge 2$. Let $\cS$ be an operator
 (essentially due to \cite{Zhao}) defined as   
\[ 
 \cS g_q(s,s-\n,a):=g_q(s,s-\n,a)+(1-q)g_q(s-1,s-\n-1,a).
\]
Since it can be easily checked that $g_q(s,s-\n-1,a)=\cS g_q(s,s-\n,a)$,
 we have inductively  
\begin{equation}\label{for:shift} 
 g_q(s,s-\n,a)=\cS^{\n-1}g_q(s,s-1,a)=\sum^{\n-1}_{r=0}\binom{\n-1}{r}(1-q)^rg_q(s-r,s-r-1,a). 
\end{equation}
Letting $q\uparrow 1$, we have $\limup g_q(s,s-\n,a)=\limup
 g_q(s,s-1,a)=\z(s,a)$. Thus the claim follows. 

 We next show the necessity. Suppose that $\limup g_q(s,\varphi(s),a)$
 exists and $\limup g_q(s,\varphi(s),a)=\z(s,a)$ for all $s\ne 1$. By
 Lemma\,\ref{lem:C123} and Proposition\,\ref{prop:q-EM}, the limit
 $(1-q)^{s-1}\frac{(1-q)}{\log{q}}b_{q^a}(\varphi(s),-s+1)$ for $q
 \uparrow 1$ does exist for all $s\ne 1$. Assume $\Re(s)<1$. Since
 $\limup(1-q)^{s-1}$ diverges, in particular, it is necessary to have
 $\limup b_{q^a}(\varphi(s),-s+1)=0$. Since $\limup
 b_{q^k}(\a,\b)=B(\a,\b)$ for all $\a, \b\in\bC$ with $B(\a,\b)$ being
 the beta function, we see that
 $B(\varphi(s),-s+1)=\G(\varphi(s))\G(-s+1)/\G(\varphi(s)-s+1)=0$. Since
 the poles of the gamma function $\G(s)$ are only at the non-positive
 integers, it follows that $\varphi(s)-s+1 \in \bZ_{\le 0}$. This shows
 that $\varphi(s)=s-\n$ for some $\n \in \bN$ for $\Re(s)<1$. Since
 $\varphi(s)$ is meromorphic in $\bC$, $\varphi(s)=s-\n$ for all
 $s\in\bC$. Hence the assertion follows. 
\end{proof}

By the discussion above, we may have infinitely many true $q$-analogues
 of the Hurwitz zeta function which are {\it not} of the form
 $f_q(s,\varphi(s),\c)$ for some $\varphi(s)$. For instance, as a
 remark, we have the following corollary.

\begin{cor}
\label{cor:tsumura}
 For a positive integer $\m$, let
\begin{align*}
 \z^{\m}_q(s,a):
&=\frac{(\m-1)!}{(1-s)_{\m}}\frac{(1-q)^s}{\log{q}}+g_q(s,\m,a)\\
&=\frac{(-1)^{\m}(\m-1)!}{(s-1)\cdots(s-\m)}\frac{(1-q)^s}{\log{q}}+\sum^{\infty}_{n=0}\frac{q^{\m(n+a)}}{{[n+a]_q}^s}.
\end{align*}
 Then the function $\z^{\m}_q(s,a)$ is a meromorphic function on $\bC$
 with simple poles at $s=1,2,\ldots,\m$, and gives a true $q$-analogue of
 $\z(s,a)$; $\limup \z^{\m}_q(s,a)=\z(s,a)$ holds for all $s\in\bC$
 $(s\ne 1,2,\ldots,\m)$.    
\end{cor}
\begin{proof}
 The first assertion is clear from the definition of $\z^{\m}_q(s,a)$.
 Further, by Proposition\,\ref{prop:q-EM} and Lemma\,\ref{lem:C123}, it
 is easy to observe the function $G_q(s,t,a)$ defined as  
\begin{equation}
\label{for:more true}
 G_q(s,t,a):=g_q(s,t,a)+\frac{(1-q)^s}{\log{q}}b_{q^a}(t,-s+1)-h_q(s,t,a)
\end{equation}
 gives a true $q$-analogue of $\z(s,a)$ if
 $\limup h_q(s,t,a)=\frac{1}{s-1}a^{-s+1}$. Since 
\[
 b_{q^a}(\m,-s+1)
=\sum^{\m}_{l=1}\frac{(-1)^{l-1}(1-\m)_{l-1}}{(1-s)_l}q^{a(\m-l)}(1-q^a)^{-s+l}
+\frac{(-1)^{\m}(\m-1)!}{(s-1)\cdots(s-\m)},
\] 
 we have $\z^{\m}_q(s,a)=G_q(s,\m,a)$, where
\[
 h_q(s,t,z)
=\sum^{\m}_{l=1}\frac{(-1)^{l-1}(1-\m)_{l-1}}{(1-s)_l}q^{a(\m-l)}\Bigl(\frac{1-q^a}{1-q}\Bigr)^{-s+l}\frac{(1-q)^l}{\log{q}}.
\]
 Hence the claim follows immediately from the fact 
 $\limup h_q(s,t,a)=\frac{1}{s-1}a^{-s+1}$.  
\end{proof}

\begin{remark}
 The function $\z^1_q(s,a)$ gives $\widetilde{\z}_q(s,a)$ in
 \cite{Tsumura2001}. 
\end{remark}


\begin{thm}\label{thm:DL}
 Let $\c$ be a Dirichlet character modulo $N$. Let $t=\varphi(s)$ be a
 meromorphic function on $\bC$. Then the formula   
\[
 \limup f_q(s,\varphi(s),\c)=L(s,\c) \quad (s \in \bC)
\]
 holds if and only if the function $\varphi(s)$ can be written as

\  $(i)$\ $\varphi(s)=s-\n$ $(\n\in\bN)$ if $\c$ is the principal
 character.  

 $(ii)$\ $\varphi(s)=s-\n$ $(\n\in\bN)$ or $\varphi(s)=\m$ $(\m\in\bN)$
 if $\c$ is not the principal character. 
\end{thm}  
\begin{proof}
 The sufficiency follows immediately from Theorem\,\ref{thm:hz} and 
 Corollary\,\ref{cor:tsumura} by the formulas    
\begin{equation}\label{for:h->L}
 f_q(s,t,\c)=\frac{1}{[N]^s_q}\sum^{N}_{k=1}\c(k)g_{q^N}(s,t,\frac{k}{N})
               \quad \textrm{and} \quad 
     L(s,\c)=\frac{1}{N^s}\sum^{N}_{k=1}\c(k)\z(s,\frac{k}{N}). 
\end{equation}
 Therefore, it suffices to show the necessity. The claim $(i)$ is
 obvious from Theorem\,\ref{thm:hz}. Thus we assume $\c$ is not the
 principal. Suppose $\limup f_q(s,t,\c)=L(s,\c)$ for all
 $s\in\bC$. Then, by Proposition\,\ref{prop:q-EM}, \eqref{for:EM-hz} and
 \eqref{for:h->L}, it holds that 
\begin{equation}\label{for:lim1}
  -\limup\frac{(1-q)^s}{N\log{q}}\sum^N_{k=1}\c(k)b_{q^k}(t,-s+1)
=\frac{1}{N}\frac{1}{s-1}\sum^N_{k=1}\c(k)k^{1-s}
 \qquad (s\in\bC). 
\end{equation}
Similarly to \eqref{for:ac-beta}, we have for any integer $M\ge 2$
\[
 b_q(\a,\b)=\sum^{M-1}_{l=1}\frac{(-1)^l(1-\a)_{l-1}}{(\b)_l}q^{\a-l}(1-q)^{\b+l-1}+\frac{(-1)^{M-1}(1-\a)_{M-1}}{(\b)_{M-1}}b_q(\a-M+1,\b+M-1)  
\]
when $\Re(\a)>M-1$. Thus we have
\begin{align}\label{for:b_q(t,-s+1)}
  b_{q^k}(t,-s+1)
=\frac{1}{s-1}q^{k(t-1)}(1-q^k)^{-s+1}+\sum^{M-1}_{l=2}\frac{(-1)^l(1-t)_{l-1}}{(-s+1)_l}q^{k(t-l)}(1-q^k)^{-s+l}\\
+\frac{(-1)^{M-1}(1-t)_{M-1}}{(-s+1)_{M-1}}b_{q^k}(t-M+1,-s+M).\nonumber
\end{align}
 Hence, using \eqref{for:b_q(t,-s+1)}, we find the formula \eqref{for:lim1} is
 equivalent to  
\begin{equation}\label{for:lim2}
 (1-t)_{M-1}\limup (1-q)^{s-1}\sum^N_{k=1}\c(k)b_{q^k}(t-M+1,-s+M)=0
 \qquad (s\in\bC)
\end{equation}
 for some $M\ge 2$. We divide the proof into the following two cases;

\begin{enumerate}
 \item The case where $t=\varphi(s)=\m\in\bN$: 

 Take $M$ as $M=t+1$. By \eqref{for:ac-beta}, we have  
\begin{multline*}
 b_{q^k}(t-M+1,-s+M)
=\frac{1}{t-(M-1)}q^{k(t-M+1)}(1-q^k)^{-s+M-1}\\
-\frac{s-(M-1)}{t-(M-1)}b_{q^k}(t-M+2,-s+M-1).\qquad\qquad
\end{multline*}
 Note that 
\[
 b_{q^k}(1,-s+M-1)=\frac{1}{s-(M-1)}\bigl((1-q^k)^{-s+M-1}-1\bigr).
\]
 Hence we have $\bigl(t-(M-1)\bigr)b_{q^k}(t-M+1,-s+M)|_{t=M-1}=1$.
 Since $\sum^N_{k=1}\c(k)=0$, the claim follows.

 \item Otherwise:

 It is clear that the condition \eqref{for:lim2} is equivalent to 
\[
 \sum^N_{k=1}\c(k)k^{1-s}\limup (1-q^k)^{s-1}b_{q^k}(t-M+1,-s+M)=0
\]
 By \eqref{for:ac-beta}, we see that for any integer $M'\ge 2$ 
\begin{align}
\label{for:M'}
& b_q(t-M+1,-s+M)=B(t-M+1,-s+M)\\
&\ \ \ \ \ +\sum^{M'-1}_{l=1}\frac{(-1)^{l-1}(1+s-M)_{l-1}}{(t-M+1)_l}q^{t-M+l}(1-q)^{-s+M-l}\nonumber\\
&\ \ \ \ \  -\frac{(-1)^{M'-1}(1+s-M)_{M'-1}}{(t-M+1)_{M'-1}}\int^1_qu^{t-M+M'-1}(1-u)^{-s+M-M'}du.\nonumber
\end{align}
 Put, $M'=M-1$ in \eqref{for:M'}. Then we obtain  
\begin{multline}
\label{for:lim3}
  \sum^{N}_{k=1}\c(k)k^{1-s}(1-q^k)^{s-1}b_{q^k}(t-M+1,-s+M)\\
=B(t-M+1,-s+M)\sum^N_{k=1}\c(k)k^{1-s}(1-q^k)^{s-1}\qquad\qquad\qquad\qquad\qquad\\
+\sum^N_{k=1}\c(k)k^{1-s}\sum^{M-2}_{l=1}\frac{(-1)^{l-1}(1+s-M)_{l-1}}{(t-M+1)_l}q^{k(t-M+l)}(1-q^k)^{M-l-1}\\
-\frac{(-1)^{M-2}(1+s-M)_{M-2}}{(t-M+1)_{M-2}}\sum^N_{k=1}\c(k)k^{1-s}\frac{(1-q^k)^s}{1-q^k}\int^1_{q^k}u^{t-2}(1-u)^{-s+1}du.
\end{multline}
 Note that we have
\[
 \frac{1}{1-q^k}\int^1_{q^k}u^{t-2}(1-u)^{-s+1}du=p^{k(t-2)}(1-p^k)^{-s+1}
\]
 for some $q<p<1$ by the mean-value theorem. Hence we have 
 $\limup \frac{(1-q^k)^s}{1-q^k}\int^1_{q^k}u^{t-2}(1-u)^{-s+1}du=0$ for
 all $s\in\bC$. Therefore, \eqref{for:lim2} is equivalent to 
\[
 B(t-M+1,-s+M)\times\biggl\{\limup (1-q)^{s-1}
 \sum^N_{k=1}\c(k)k^{1-s}\Bigl(\frac{1-q^k}{1-q}\Bigr)^{s-1}\biggr\}=0
 \qquad(s\in\bC). 
\]
 Note that there exists always an integer $r\ge 0$ such that
 $\sum^f_{k=1}\c(k)(k-1)^r\ne 0$. Assume $\Re(s)<1-r$. Since $\limup
 (1-q)^{s-1}\sum^N_{k=1}\c(k)k^{1-s}\bigl(\frac{1-q^k}{1-q}\bigr)^{s-1}$
 diverges for such $s$ by the following lemma, we conclude that 
 $B(t-M+1,-s+M)=0$. Hence the claim follows from the same discussion as
 in the proof of Theorem\,\ref{thm:hz}.
\end{enumerate}        
\end{proof}


\begin{lem}
 We have
\begin{equation}\label{for:asym}
 \Bigl(k^{-1}\frac{1-q^k}{1-q}\Bigr)^{s-1}=\sum^r_{i=0}(q-1)^i\binom{s-1}{i}
 \Bigl(\sum^{k-1}_{l=1}[l]_q\Bigr)^ik^{-i}+O((1-q)^{r+1})
\end{equation}
for $r\ge0$. In particular,  $\limup \Bigl|
 (1-q)^{s-1}\sum^N_{k=1}\c(k)k^{1-s}\bigl(\frac{1-q^k}{1-q}\bigr)^{s-1}\Bigr|=\infty$
 for $\Re(s)<1-r$ provided $\sum^N_{k=1}\c(k)(k-1)^r\ne 0$. 
\end{lem}
\begin{proof}
 Since
 $k^{-1}\frac{1-q^k}{1-q}=1+([1]_q+[2]_q+\cdots+[k-1]_q)(q-1)k^{-1}$,
 \eqref{for:asym} follows from the binomial expansion. The rest of the
 assertion follows immediately from \eqref{for:asym}. 
\end{proof}


\begin{remark}
 Let $f(z)$ be an arbitrary function of the form
 $f(z)=z^{-s+1}P(z)+z^{-1}Q(z^{-1})\d(\c)$ with $f(1)=1$, $P(z)$ being 
 a polynomial and $Q(z)$ being a holomorphic function at $z=0$. Here
 $\d (\c)=1$ if $\c$ is not the principal character and $\d (\c)=0$
 otherwise. Then, it is clear from Theorem\,\ref{thm:DL} that
 $\limup L_f(s,\c)=L(s,\c)$ where
 $L_f(s,\c):=\sum^{\infty}_{n=1}\c(n)f(q^{-n})[n]_q^{-s}$. Note that,
 however, $P$ can not be a (infinite) Taylor series because the series
 $L_f(s,\c)$ does not converge.    
\end{remark}

Special values of $L\nn_q(s,\c)$ are obtained by the following
 expression \eqref{for:L_q} of $L\nn_q(s,\c)$, which gives also the
 information concerning the locations of the poles (and a meromorphic
 continuation again) of $L\nn_q(s,\c)$.     
 

\begin{prop}\label{prop:specialvalue}
 (i)\ Let $\c$ be a Dirichlet character. Then $L\nn_q(s,\c)$ can be
 written as 
\begin{equation}\label{for:L_q}
 L\nn_q(s,\c)=(1-q)^s\sum^{\infty}_{r=0}\binom{s+r-1}{r}g_{\c}(q^{s-\n+r}),
\end{equation}
 where $g_{\c}(q)=\sum^N_{k=1}\c(k)q^k/(1-q^N)$. In particular, if $\c$
 is not the principal character, then $L\nn_q(s,\c)$ is holomorphic at
 $s=1,2,\ldots,\n$, and if $\c$ is the principal character, then
 $L\nn_q(s,\c)$ has simple poles at these points. Especially,
 $\z\nn_q(s)$ has simple poles at these points.\\
(ii)\ Let $m$ be a non-negative integer. Then we have
\begin{equation}\label{for:specialvalue}
 L\nn_q(-m,\c)=(1-q)^{-m}\left\{\sum^m_{r=0}(-1)^r
\binom{m}{r}g_{\c}(q^{-m+r-\n})+\frac{(-1)^{m+1}m!(\n-1)!}{(m+\n)!\log{q}}B_{0,\c}\right\},
\end{equation}
where $B_{n,\c}$ is the generalized Bernoulli number defined via
 $\sum^{N}_{k=1}\frac{\c(k)te^{kt}}{e^{Nt}-1}=\sum^{\infty}_{n=0}B_{n,\c}\frac{t^n}{n!}$.  
\end{prop}
\begin{proof}
 Using the binomial theorem
 $(1-x)^{-s}=\sum^{\infty}_{r=0}\binom{s+r-1}{r}x^r$, we can show 
\begin{equation}\label{for:f_q}
 f_q(s,t,\c) =(1-q)^s\sum^{\infty}_{r=0}\binom{s+r-1}{r}g_{\c}(q^{t+r})
\end{equation}
 by the same way as in \cite{KanekoKurokawaWakayama2003}. Indeed, the
 change of the order of two summations can be justified because the
 series converge absolutely. We notice here that
\begin{equation}\label{for:g_c}
 g_{\c}(q^{t+r})=-\sum^{\infty}_{n=-1}B_{n+1,\c}\frac{(\log{q})^n}{(n+1)!}(t+r)^n. 
\end{equation}
 By the help of \eqref{for:g_c}, the equation
 \eqref{for:f_q} gives the rest of the first assertion. The formula
 \eqref{for:specialvalue} follows from \eqref{for:L_q}.    
\end{proof}

From \eqref{for:specialvalue} and \eqref{for:g_c} we have the well-known
formula $L(-m,\c)=-B_{m+1,\c}/(m+1)$ for a non-negative integer $m$ as
the classical limit $q \uparrow 1$. We next consider the limit $q
\downarrow 0$ of $L\nn_q(s,\c)$ which we call a {\it crystal limit}. We
then introduce functions $L\nn_0(s,\c)$ and $\z\nn_0(s)$ by the
point-wise limits;    
\[
 L\nn_0(s,\c):=\limdo L\nn_q(s,\c) \quad \textrm{and} \quad
 \z\nn_0(s):=\limdo \z\nn_q(s). 
\]
We call the function $\z\nn_0(s)$ (resp. $L\nn_0(s,\c)$) a {\it crystal
Riemann zeta} (resp. $L$-) {\it function of type $\n$}. Note that, in
particular, since $\z\nn_q(s)$ has simple poles at $s=1,2,\ldots,\n$,
$\z\nn_0(s)$ can not be defined at these points. The following 
proposition is easily obtained from Proposition\,\ref{prop:specialvalue}.      


\begin{prop}\label{prop:crystal}
 Let $m$ be a non-negative integer. Then we have 
\begin{align}\label{for:crystal}
  L\nn_0(-m,\c)&=0 \quad (\c\ne{\bf 1}) \quad \textrm{and} \quad 
  \z\nn_0(-m)=\begin{cases}
 -1 & (m=0),\\
  0 & (m\ne 0).
\end{cases}
\end{align}
Further, let $D\nn_{0}:=\{s\in\bC\,|\,\Re(s)\not\in \bZ_{\le
 \n}\}\cup \{0,-1,-2,\ldots\}$. Then for $s \in D\nn_{0}$, it holds
 that  
\begin{align*}
\ \  L\nn_0(s,\c)&=0 \qquad (\c\ne {\bf 1}),\\
  \z\nn_0(s)&=
\Biggl\{
\begin{array}{ccl}
  0             & \textrm{if} &\     \Re(s)>\n,    \\ 
  -(s+1)_{m}/m! & \textrm{if} &\  \n-m-1<\Re(s)<\n-m \quad (m=0,1,2,\ldots).
\end{array}
\Biggr.
\end{align*}
Note that $\limdo L\nn_q(s,\c)$ does not exist if $s\notin D\nn_{0}$.\qed
\end{prop}

As a corollary of Proposition\,\ref{prop:crystal}, we now determine the
zeros of $\z\nn_0(s)$. This can be regarded as analogues of the Riemann
hypothesis for the crystal Riemann zeta functions. 


\begin{cor}\label{cor:crystal-RH}
 If $\z\nn_0(s)=0$, $\Re(s) \le \n$ and $s \ne 1,2,\ldots,\n$, then
 $s \in \bZ_{<0}$. \qed 
\end{cor}

By this corollary, it seems that the crystal zeta functions have only
``trivial'' zeros (See the figures in \S 3). In general, however, no
reason can judge which zeros are the trivial or the non-trivial zeros of
$\z\nn_q(s)$ a priori because we do not have a functional equation of
$\z\nn_q(s)$. In the next section, we therefore study numerically the
zeros of $\z\nn_q(s)$ as $q$-trajectories of the trivial and non-trivial
zeros of $\z(s)$, respectively, when $q$ approaches $0$ decreasing from
$1$.     


\begin{remark}
 Let $\c$ be a even character, that is, the equation
 $\sum^N_{k=1}k\c(k)=0$ holds. Then $L(1,\c)$ is expressed as $L'(0,\c)$
 by the functional equation of $L(s,\c)$. In our case,  
\begin{align*}
 L\nn_q{'}(0,\c)
&=\sum^N_{k=1}\c(k)q^{-k\n}\left\{\frac{\log{q}(k+(N-k)q^{-N\n})}{(1-q^{-N\n})^2}+\frac{\log{(1-q)}}{1-q^{-N\n}}+\sum_{1
\le r \ne \n}\frac{1}{r}\frac{q^{kr}}{1-q^{N(r-\n)}}\right\},\\ 
 L\nn_q(1,\c)
&=(1-q)\sum^N_{k=1}\c(k)q^{-k\n}\left\{\sum_{1 \le r \ne \n}\frac{q^{kr}}{1-q^{N(r-\n)}}\right\}.
\end{align*}
Hence it is difficult to express $L\nn_q(1,\c)$ in terms of
 $L\nn_q{'}(0,\c)$. This partially suggests that one may not expect the
 presence of a functional equation of $L\nn_q(s,\c)$.   
\end{remark}


\begin{remark}
 Let $\c$ be a Dirichlet character modulo $N$ but the principal
 character. Then the Dirichlet class number formula is given by 
\begin{equation}\label{for:L1}
 L(1,\c)=-\frac{1}{N}\sum^N_{k=1}\c(k)\frac{\G'}{\G}(\frac{k}{N}).
\end{equation}
In our case, similarly, $L\nn_q(1,\c)$ can be expressed in terms of
 $\G_q(s)$ as  
\begin{equation}\label{for:q-L1}
 L\nn_q(1,\c)=(1-q)\sum^{\n-1}_{r=1}g_{\c}(q^{-\n+r})+\frac{1-q}{N\log{q}}\sum^N_{k=1}\c(k)\frac{\G'_{q^N}}{\G_{q^N}}(\frac{k}{N}).
\end{equation}
Here $\G_q(s)$ is the Jackson $q$-gamma function defined by
 $\G_q(s):=\frac{(q;q)_{\infty}}{(q^s;q)_{\infty}}(1-q)^{1-s}$, where 
 $(a;q)_n=\prod^{n-1}_{j=0}(1-aq^j)$ (see,
 \cite{AndrewsAskeyRoy1999}). This is shown by comparing $L\nn_q(1,\c)$
 with the logarithmic derivative of $\G_q(s)$. Letting $q \uparrow 1$ in
 \eqref{for:q-L1}, we obtain \eqref{for:L1} by
 Theorem\,\ref{thm:DL}. Though a certain quantum analogue of the
 Dirichlet class number formula was studied in
 \cite{KurokawaWakayama2002}, we remark that the $q$-analogue of the
 Hurwitz zeta function defined there is not the same one discussed
 here. See also \cite{KurokawaWakayama2003}.       
\end{remark}

%
\section{Behavior of zeros of $\z\nn_q(s)$}
%

In this section, we study the zeros of $\z\nn_q(s)$ numerically
with Maple $8$ \cite{maple} and give some conjectures concerning the
$q$-trajectory of the zeros of $\z(s)$.

We first notice the 


\begin{prop}
 For all $q\in (0,1]$, $\z\nn_q(s)\ne 0$ if $\Re(s)\ge 2\n$.
\end{prop}
\begin{proof}
 Since $\z\nn_q(s)=q^{s-\n}+\sum^{\infty}_{n=2}q^{n(s-\n)}[n]^{-s}_q$,
 it is sufficient to show
 $\sum^{\infty}_{n=2}|q^{n(s-\n)}[n]^{-s}_q|<|q^{s-\n}|$. To
 see this, since $[n]_q=1+q+\cdots+q^{n-1}\ge n(1\cdot q\cdot\cdots\cdot
 q^{n-1})^{\frac{1}{n}}=nq^{\frac{n-1}{2}}$, it is enough to verify
 $\sum^{\infty}_{n=2}q^{(n-1)(\frac{\s}{2}-\n)}n^{-\s}<1$, where
 $\s=\Re(s)$. This is actually true because    
\[
 \sum^{\infty}_{n=2}q^{(n-1)(\frac{\s}{2}-\n)}n^{-\s}\le\sum^{\infty}_{n=2}n^{-\s}<
 \int^{\infty}_{1}x^{-\s}dx <\frac{1}{\s-1}\le1 
\]
for $\s\ge 2\n$. This shows the assertion.
\end{proof}

We begin with a study of an approximate formula of $\z\nn_q(s)$. Put
$f(x)=q^{xt}(1-q^x)^{-s}$. Then by the Euler-Maclaurin formula
\eqref{for:EM} again, we have 
\begin{align}\label{for:nc-f_q}
 f_q(s,t)
 &=(1-q)^s\biggl\{\sum^{N}_{m=1}q^{mt}(1-q^m)^{-s}-\frac{1}{2}q^{Nt}(1-q^N)^{-s}-
   \frac{1}{\log{q}}b_{q^N}(t,-s+1) \biggr.\\
 &\ \ \ \biggl.-\sum^{n}_{k=1}\frac{B_{2k}}{(2k)!}f^{(2k-1)}(N) \biggr\} -
  \frac{(1-q)^s}{(2n)!}\int^{\infty}_{N}\widetilde{B}_{2n}(x)f^{(2n)}(x)dx.\nonumber  
\end{align}
Next, we try to evaluate the integral in terms of the incomplete beta
functions. For a non-negative integer $j$, we define the function
$h_j(x)$ by $q^{(t+j)x}(1-q^x)^{-s-j}$. Then $h^{(n)}_0(x)(=f^{(n)}(x))$
is expressed by a linear combination of $h_j(x)$'s with coefficients   
\[
 a^{(n)}_j := (s)_j\sum_{i_0+i_1+\cdots+i_j \atop
 =n-j}t^{i_0}(t+1)^{i_1}\cdots(t+j)^{i_j}. 
\]


\begin{lem}\label{lem:h}
For any $n \ge 0$, we have 
\[
h^{(n)}_0(x)=(\log{q})^n\sum^n_{j=0}a^{(n)}_jh_j(x).
\]
\end{lem}
\begin{proof}
 The diagram (based on Leibniz rule) below describes the rule how the
  coefficient of $h_j(x)$ in
  $h^{(n)}_0(x)=\left(q^{tx}(1-q^x)^{-s}\right)^{(n)}$ can be obtained.   
\[
 \xymatrix{
    & & & h_0 \ar[dl]_t \ar[dr]^s & & &\\
    & & h_0 \ar[dl]_t \ar[dr]^s & & h_1 \ar[dl]_{t+1} \ar[dr]^{s+1} &\\
    & h_0 \ar[dl]_t \ar[dr]^s & & h_1 \ar[dl]_{t+1} \ar[dr]^{s+1} & & h_2
 \ar[dl]_{t+2} \ar[dr]^{s+2} &\\ 
    h_0 & & h_1 & & h_2 & &  h_3 &  }
\]
Hence the claim follows immediately.
\end{proof}

By \eqref{for:Fourier}, \eqref{for:ac-beta} again and
Lemma\,\ref{lem:h}, for any integer $M \ge 2$ and $\Re(t)>1-M$, the last 
term of \eqref{for:nc-f_q} containing the integral can be evaluated as  
\begin{align*}
& (1-q)^s (\log{q})^{2n-1} \Biggr. \sum_{l_0 \le l \le l_1 \atop
 l\ne0}\sum^{2n}_{j=0}\sum^{M-1}_{k=1}\frac{a^{(2n)}_{j}}{(2\pi
 il)^{2n}}\frac{(-1)^{k-1}(s+j)_{k-1}}{(t+j+\d
 l)_k}\frac{q^{N(t+j+k-1)}}{(1-q^N)^{s+j-1+k}}\Biggl.\\
&+ R(s,t,q,N,M,n,l_0,l_1). 
\end{align*}
Here $l_0$ and $l_1$ are integers satisfying $l_0 < l_1$, and
$R(s,t,q,N,M,n,l_0,l_1)$ is equal to 
\begin{align}\label{for:reminder}
&(\log{q})^{2n-1}(1-q)^s\Biggl\{\sum_{l<l_0, l>l_1 \atop
 l\ne0}\sum^{2n}_{j=0}\sum^{M-1}_{k=1}\frac{a^{(2n)}_{j}}{(2\pi il)^{2n}}\frac{(-1)^{k}(s+j)_{k-1}}{(t+j+\d
 l)_k}\frac{q^{N(t+j+k-1)}}{(1-q^N)^{s+j-1+k}}\Biggl.\\
&\ \ \ +\Biggl.\sum_{l\in\bZ\backslash\{0\}}\sum^{2n}_{j=0}\frac{a^{(2n)}_{j}}{(2\pi il)^{2n}}\frac{(-1)^{M}(s+j)_{M-1}}{(t+j+\d l)_{M-1}}b_{q^N}(t+j+\d l+M-1,-s-j-M+2)\Biggl\}\nonumber.
\end{align}
Substituting $s-\n$ for $t$ in $f_q(s,t)$, we obtain the following
expression which allows us to calculate the zeros of $\z\nn_q(s)$ 
numerically.    


\begin{prop}
 For integers $N \ge 1, M \ge 2, n \ge 1$ and $l_0, l_1$ satisfying
 $l_0<l_1$, we have for $\Re(s)>\n+1-M$   
\begin{align*}
&\ \ \z\nn_q(s)=(1-q)^s\Biggl\{\sum^{N}_{m=1}q^{m(s-\n)}(1-q^m)^{-s}-\frac{1}{2}q^{N(s-\n)}(1-q^N)^{-s}\\
&\ +\frac{1}{\log{q}}\sum^{\n-1}_{r=0}\frac{(-\n+1)_r}{(-s+1)_{r+1}}q^{N(s-\n)}(1-q^N)^{-s+1+r}-\sum^n_{k=1}\sum^{2k-1}_{j=0}\frac{B_{2k}}{(2k)!}(\log{q})^{2k-1}a^{(2k-1)}_jq^{N(s-\n+j)}(1-q^N)^{-s-j}\\
&\ \ \Biggl.+(\log{q})^{2n-1}\sum_{l_0 \le l \le l_1 \atop
 l\ne0}\sum^{2n}_{j=0}\sum^{M-1}_{k=1}\frac{a^{(2n)}_{j}}{(2\pi
 il)^{2n}}\frac{(-1)^{k}(s+j)_{k-1}}{(s-\n+j+\d
 l)_k}\frac{q^{N(s-\n-1+j+k)}}{(1-q^N)^{s+j-1+k}}\Biggl\}\\
&\ \ \ +R(s,s-\n,q,N,M,n,l_0,l_1).
\end{align*}
Moreover, the absolute value of $R(s,s-\n,q,N,M,n,l_0,l_1)$ is bounded by
\begin{align*}
&\frac{|\log{q}|^{2n-1}(1-q)^{\Re(s)}}{(2\pi)^{2n}}\Biggl\{\sum_{l<l_0,
 l>l_1 \atop
 l\ne0}\sum^{2n}_{j=0}\sum^{M-1}_{k=1}\frac{|a^{(2n)}_{j}|}{l^{2n}}\frac{|(s+j)_{k-1}|}{|(s-\n+j)_k|}\frac{q^{N(\Re(s)-\n-1+j+k)}}{(1-q^N)^{\Re(s)+j-1+k}}\Biggl.\\ 
&\ \ \
 -\sum_{l\in\bZ\backslash\{0\}}\sum^{2n}_{j=0}\frac{|a^{(2n)}_{j}||(s+j)_{M-1}|}{l^{2n}|(\Re(s)-\n+j)_{M-1}|}\sum^{\n-1}_{r=0}\frac{(-\n+1)_r}{(-\Re(s)-j-M+2)_{r+1}}\frac{q^{\Re(s)+j+M-\n-1}}{(1-q)^{\Re(s)+j+M-2-r}}\Biggl\}. 
\end{align*} 
\end{prop}
\begin{proof}
 Let $t=s-\n$. Then the first assertion follows immediately from
 $(3.1)$. The second assertion follows from \eqref{for:reminder} and 
 \eqref{for:special-beta}. This completes the proof. 
\end{proof}

 If $s$ is real, then $\z\nn_q(s)$ is also real by the definition of
 $\z\nn_q(s)$. We actually compute approximate values of the zeros of
 $\z\nn_q(s)$ by searching the point where the sign of $\z\nn_q(s)$
 changes. If $s$ is complex, we can not, however, determine the place of
 the zeros of $\z\nn_q(s)$ as we can not expect an existence of a
 functional equation of $\z\nn_q(s)$. Thus we try to seek the points
 $s_0$ which are expected as the zeros of $\z\nn_q(s)$ by observing
 respectively the sign of the real part and the imaginary part of
 $\z\nn_q(s)$. Indeed, if $\z\nn_q(s_0)=0$, then the signs of both 
 $\Re(\z\nn_q(s))$ and $\Im(\z\nn_q(s))$ will change in a neighborhood
 of the point $s_0$, and consequently, $|\z\nn_q(s)|$ is very small in
 the neighborhood of $s_0$. We find a certain neighborhood in which the
 both of the signs simultaneously change.  
 
Following the strategy mentioned above, we examine the sign of
$\z\nn_q(s)$ when $s$ is real, and the sign of $\Re(\z\nn_q(s))$ and of
$\Im(\z\nn_q(s))$ when $s$ is complex. In these numerical calculations,
the parameter $q$ moves from $0.99$ to $0.01$ by $0.01$, and
additionally $3$ points $0.001$, $0.0001$, and $0.00001$. The parameter
$s$ moves from the real zero point at $q=1$ by $0.001$ on the points of
the real axis if $s$ is real, and on the points of the lattice parallel
to the real and imaginary axis if $s$ is complex. Furthermore we may
carefully select parameters $q, N, M, n, l_0, l_1$ satisfying
$|R(s,s-\n,q,N,M,n,l_0,l_1)|<10^{-5}$ on each point $s$. In the
following figures, the point $\times$ represents the zeros of
$\z(s)$. When $q=0.99, \ldots, 0.01$, we plot approximate (if $s$ is
real) or expecting points (if $s$ is complex) by $\bullet$, and when $q$
is the other $3$ points, by $\star$. Further,
\textcolor{red}{$\bullet$}, \textcolor{green}{$\bullet$} and
\textcolor{blue}{$\bullet$} denote $\n=1,2$ and $3$, respectively.  

\smallbreak
\smallbreak


\noindent
{\it 1. Trivial zero's cases:}\\ 
First we observe that the behavior of the zeros of $\z\nn_q(s)$
$(q\in(0,1])$ starting from the trivial zeros of $\z(s)$. Put $s_j=-2j$, 
the $j$-th trivial zero of $\z(s)$. Let $s\nn_j(q)$ be the
$q$-trajectory of  $s_j=s\nn_j(1)$ satisfying $\z\nn_q(s\nn_j(q))=0$. We
plot $s\nn_j(q)$ approximately on the ($s,q$)-plane. Note that
$s\nn_j(q)$ is real. Figure $1$, $2$ show the $q$-trajectory of
$s\nn_1(q)$, $s\nn_2(q)$ respectively. It can be expected from these two
figures that $s\nn_1(q)$, $s\nn_2(q)$ approach $s=-1, -2$ respectively
for {\it all} $\n=1,2,3$. Moreover, by virtue of further numerical
calculations, it seems true that $s\nn_3(q)$, $s\nn_4(q)$ approach
$s=-3$, $-4$ respectively for {\it all} $\n=1,2,3$. This motivates the    

\smallbreak
\smallbreak
\noindent
\textbf{Conjecture.}\textit{\ The limit $s\nn_j(0):=\limdo
s\nn_j(q)$ exists and is given by $s\nn_j(0)=s_j/2=-j$ for all
$\n\in\bN$. \qed} 
\smallskip 
 
\begin{figure}[htbp]
 \centering
 \psfrag{q}{$q$}
 \psfrag{s}{$s$}
 \includegraphics[clip]{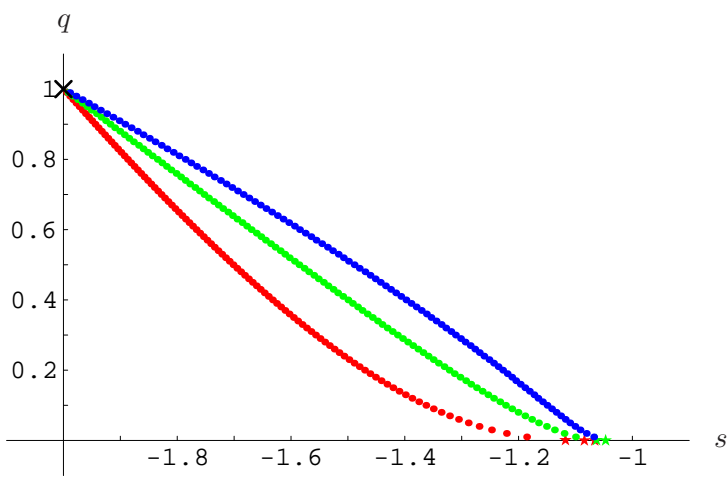}
 \caption{the trajectory of zeros of $\z\nn_q(s)$ starting from $s_1=-2$}
\end{figure}

\begin{figure}[htbp]
 \centering
 \psfrag{q}{$q$}
 \psfrag{s}{$s$}
  \includegraphics[clip]{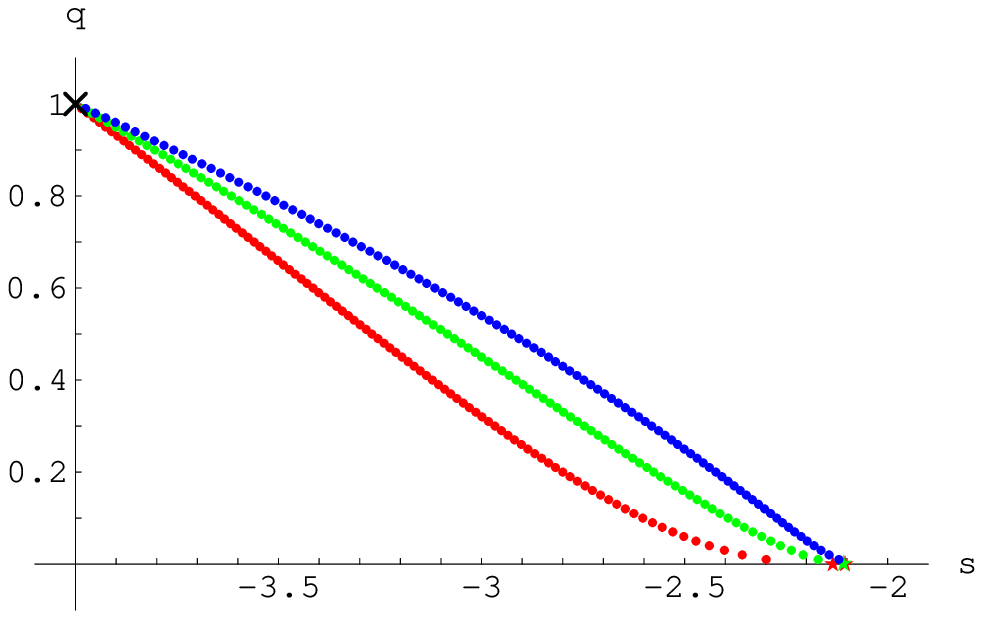}
 \caption{the trajectory of zeros of $\z\nn_q(s)$ starting from $s_2=-4$}
\end{figure}

The next proposition may partially support this conjecture.


\begin{prop}\label{prop:tangent}
 Let $z\nn(q)$ be a zero of $\z\nn_q(s)$. If $\frac{\p}{\p
 s}\z\nn_q(z\nn(q))\ne 0$ and $z\nn(q)$ is right continuous at $q=0$ and
 $z\nn(0):=\limdo z\nn(q)=-m$ for an integer $m$, then $\limdo
 \frac{d}{dq}z\nn(q)=\pm\infty$. In particular, $z\nn(q)$ is tangent to
 the line $q=0$ at $s=-m$ in $(s,q)\in\bR\times[0,1]$.   
\end{prop}
\begin{proof}
 Note that $\z\nn_q(s)$ is real for $s\in\bR$. Since $\frac{\p}{\p
 s}\z\nn_q(z\nn(q))\ne 0$, we conclude that $z\nn(q)$ is differentiable
 near $q=0$ and $\frac{d}{dq}z\nn(q)=-\frac{\p}{\p
 q}\z\nn_q(z\nn(q))\frac{\p}{\p s}\z\nn_q(z\nn(q))^{-1}$ by the implicit 
 function theorem. Since, by the expression \eqref{for:L_q},
 $\frac{\p}{\p q}\z\nn_q(s)$ and $\frac{\p}{\p s}\z\nn_q(s)$ can be
 written as   
\begin{align*}
 \frac{\p}{\p q}\z\nn_q(s)
&=-\frac{s}{1-q}\z\nn_q(s)+q^{-1}(1-q)^{-m}\\
&\ \ \ \times\Biggl\{\frac{(-1)^m m!(\n-1)!}{(\n+m)!(\log{q})^2}
+\sum^{\infty}_{r=0 \atop
 r\ne\n+m}\binom{-m+r-1}{r}\frac{(-m+r-\n)q^{-m+r-\n}}{(1-q^{-m+r-\n})^2}\Biggr\}+h_1(s),\\ 
\frac{\p}{\p s}\z\nn_q(s)
&=\log{(1-q)}\z\nn_q(s)+(1-q)^{-m}\Biggl\{\frac{(-1)^{m+1}
 m!(\n-1)!}{(\n+m)!\log{q}}\Bigl(\sum^{m+\n-1}_{k=0 \atop
 k \ne m}\frac{1}{-m+k}+\frac{1}{2}\log{q}\Bigr)\\
&\ \ \
 +\sum^m_{r=0}\binom{-m+r-1}{r}\frac{q^{-m+r-\n}}{1-q^{-m+r-\n}}\Bigl(\sum^{r-1}_{k=0}\frac{1}{-m+k}+\frac{\log{q}}{1-q^{-m+r-\n}}\Bigr)\\
&\ \ \  +\sum^{\infty}_{r=m+1 \atop r \ne m+\n}\frac{(-1)^m m!(r-m-1)!}{r!}\frac{q^{-m+r-\n}}{1-q^{-m+r-\n}}\Biggr\}+h_2(s),
\end{align*}
where $h_i(s)$ is a function satisfying $h_i(-m)=0$ $(i=1,2)$, we have   
\begin{align*}
 \limdo \frac{\p}{\p q}\z\nn_q(z\nn(q))
&=\pm\infty,\\
 \limdo \frac{\p}{\p s}\z\nn_q(z\nn(q))
&=\frac{(-1)^{m+1}m!(\n-1)!}{2(\n+m)!}-\sum^{m}_{r=0}\binom{-m+r-1}{r}\sum^{r-1}_{k=0}\frac{1}{-m+k}\\
&\ \ \ -\sum^{m+\n-1}_{r=m+1}\frac{(-1)^{m+1}m!(r-m-1)!}{r!}.
\end{align*}
 This shows the assertion. 
\end{proof}

\begin{figure}[bp]
 \centering
 \psfrag{Re}{$\Re$}
 \psfrag{Im}{$\Im$}
  \includegraphics[clip]{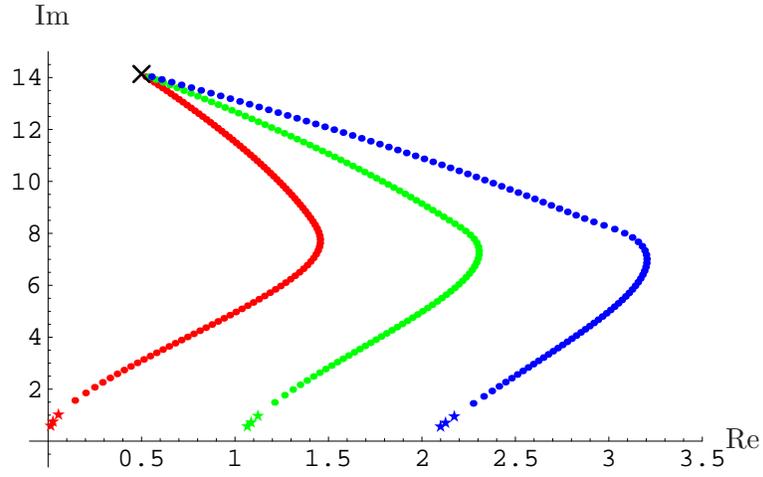}
 \caption{the trajectory of zeros of $\z\nn_q(s)$ starting from
 $\r_1=\frac{1}{2}+\sqrt{-1}\times 14.13472\ldots$}  
\end{figure}

\begin{figure}[htbp]
 \centering
 \psfrag{x}{$\Re$}
 \psfrag{y}{$\Im$}
 \psfrag{z}{$q$}
  \includegraphics[clip]{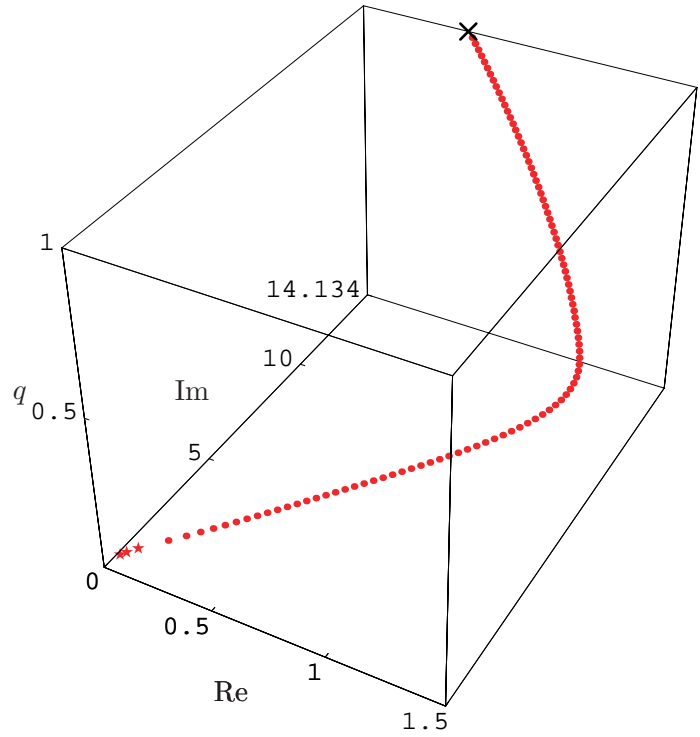}
 \caption{$3$-dimensional graph of Figure $3$ for $\n=1$}  
\end{figure}

\begin{figure}[htbp]
 \centering
 \psfrag{Re}{$\Re$}
 \psfrag{Im}{$\Im$}
  \includegraphics[clip]{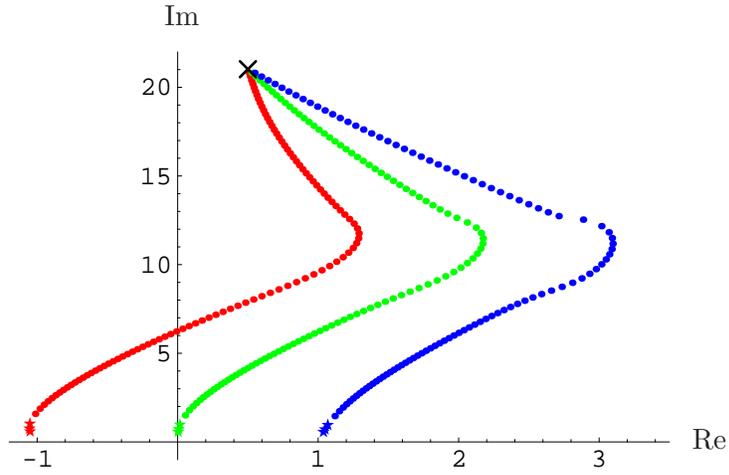}
 \caption{the trajectory of zeros of $\z\nn_q(s)$ starting from
 $\r_2=\frac{1}{2}+\sqrt{-1}\times 21.02203\ldots$}  
\end{figure}

\begin{figure}[htbp]
 \centering
 \psfrag{Re}{$\Re$}
 \psfrag{Im}{$\Im$}
  \includegraphics[clip]{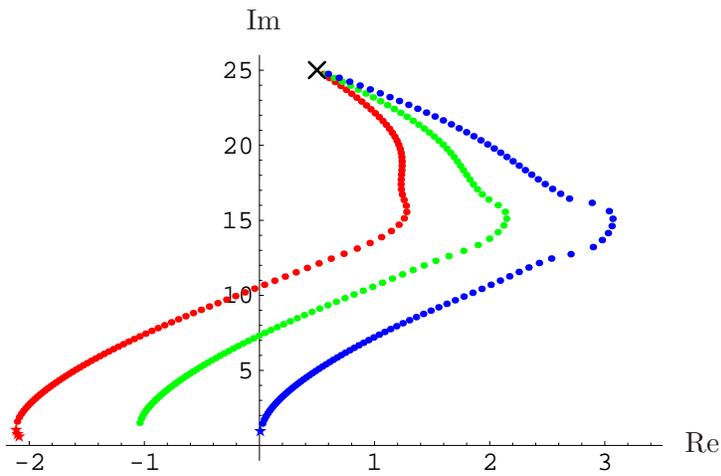}
 \caption{the trajectory of zeros of $\z\nn_q(s)$ starting from
 $\r_3=\frac{1}{2}+\sqrt{-1}\times 25.01085\ldots$}  
\end{figure}

\begin{figure}[htbp]
 \centering
 \psfrag{Re}{$\Re$}
 \psfrag{Im}{$\Im$}
  \includegraphics[clip]{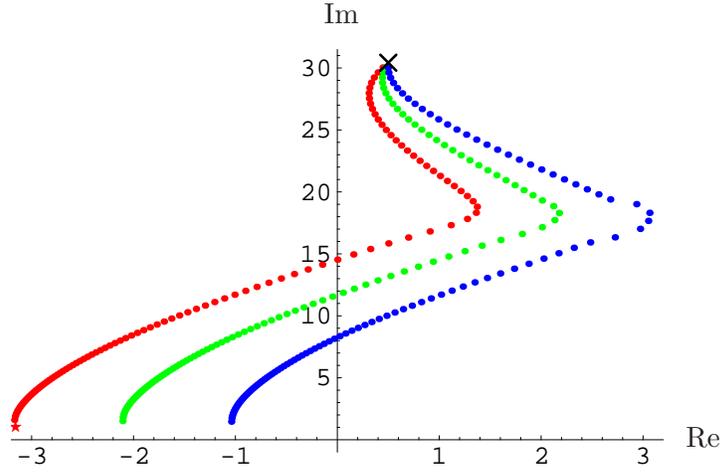}
 \caption{the trajectory of zeros of $\z\nn_q(s)$ starting from
 $\r_4=\frac{1}{2}+\sqrt{-1}\times 30.42487\ldots$}  
\end{figure}

\begin{figure}[htbp]
 \centering
 \psfrag{Re}{$\Re$}
 \psfrag{Im}{$\Im$}
  \includegraphics[clip]{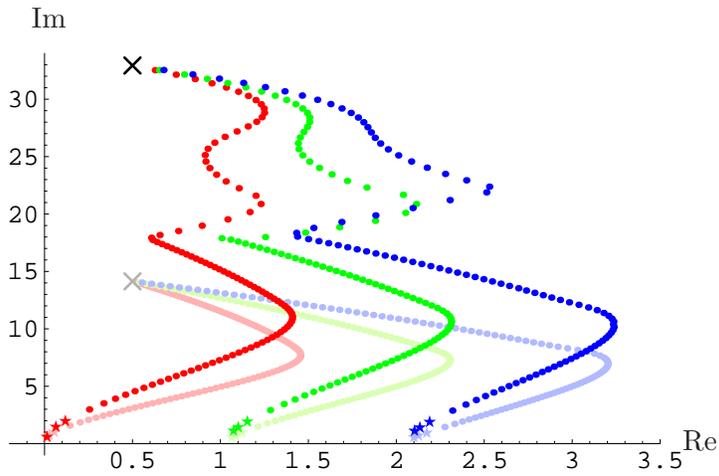}
 \caption{the trajectory of zeros of $\z\nn_q(s)$ starting from
 $\r_5=\frac{1}{2}+\sqrt{-1}\times 32.93506\ldots$
 with $\r_1$-trajectory (see Figure\,3)}  
\end{figure}


\noindent
{\it 2. Non-trivial zero's cases:}\\ 
We here investigate the behavior of the zeros of $\z\nn_q(s)$
$(q\in(0,1])$ starting from the non-trivial zeros of $\z(s)$. Let $\r_j$
$(j=1,2,\ldots)$ be the $j$-th non-trivial zero of $\z(s)$ whose
imaginary part is positive, that is, $\r_1=\frac{1}{2}+\sqrt{-1}\times
14.13472\ldots$, $\r_2=\frac{1}{2}+\sqrt{-1}\times 21.02203\ldots$, 
$\r_3=\frac{1}{2}+\sqrt{-1}\times 25.01085\ldots$, $\ldots$. Let
$\r\nn_j(q)$ be the $q$-trajectory of $\r_j=\r\nn_j(1)$ satisfying
$\z\nn_q(\r\nn_j(q))=0$. Note that, by the principle of reflection, we
have $\z\nn_q(\overline{\r\nn_j(q)})=0$. Figure $3$ (see also Figure $4$),
Figure $5$, Figure $6$, Figure $7$ and Figure $8$ correspond to the case
$\r_1$, $\r_2$, $\r_3$, $\r_4$ and $\r_5$ respectively. In these cases, we plot
$\r\nn_j(q)$ on the complex plane. $\r\nn_j(q)$ $(j=1,2,\ldots,5)$ seem
to arrive at $0$, $-1$, $-2$, $-3$ and $0$ respectively for $\n=1$ when
$q \downarrow 0$. It is easy to see that if $\n$ increases $1$, the real
part of $\r\nn_j(q)$ also increases $1$ near $q=0$.

\smallbreak
\smallbreak
\noindent
\textbf{Conjecture.}\textit{\ The limit $\r\nn_j(0):=\limdo\r\nn_j(q)$
exists and satisfies $\r\nn_j(0)\in\bZ_{\le \n-1}$. Additionally, the
relation $\r^{(\n+1)}_j(0)=\r\nn_j(0)+1$ holds for $\n\in\bN$. \qed}     


\begin{remark}
 If $\r\nn_j(q)$ is right continuous at $q=0$, then we can show that
 $\Re(\r\nn_j(0))\in\bZ_{<0}$. In fact, suppose otherwise. Then there
 exists some positive integer $N$ such that $-N<\Re(\r\nn_j(q))<-N+1$
 for $q\in [0,\e)$ for some $\e>0$. Since $\z\nn_q(\r\nn_j(q))=0$, the
 continuity of $\r\nn_j(q)$ and the existence of the limit show that
 $0=\limdo\z\nn_q(\r\nn_j(q))=\z\nn_0(\r\nn_j(0))$. This contradicts the
 results in Proposition\,\ref{prop:crystal}.  
\end{remark}


\begin{remark}\label{rem:rho5}
 Remark that the limit behavior of the $5$th zeros (Figure $8$) is
 different from the one what can be expected from the $1$st, $2$nd,
 $3$rd and $4$th zeros. See also Figures $13$, $14$, $15$, and the
 subsequent observation.   
\end{remark}


\noindent
{\it 3. Other remarks and some questions:}\\
We give another numerical calculation of $\z_q(s)=\z_q^{(1)}(s)$
at the neighborhoods of $s=0,-1$ with small $q$ from the viewpoint of the
approximations obtained by Proposition\,\ref{prop:specialvalue}:
\[
 \z_q(s) \approx 
\begin{cases}
(1-q)^{s}\Bigl(\frac{q^{s-1}}{1-q^{s-1}}+s\frac{q^{s}}{1-q^{s}}\Bigr)
 & \text{if $\Re(s) > -1$}, \\
(1-q)^{s}\Bigl(\frac{q^{s-1}}{1-q^{s-1}}+s\frac{q^{s}}{1-q^{s}} + 
\frac{s(s+1)}{2} \frac{q^{s+1}}{1-q^{s+1}} \Bigr)
 & \text{if $\Re(s) > -2$}.
\end{cases}
\]
Figure $9$ and Figure $10$ indicate the logarithm of absolute values of 
these forms on the rectangles $[-0.05,\, 0.05] \times [0,\, 1]$ and
$[-1.05,\, -0.95] \times [0,\, 1]$ with $q=2^{-64}$, respectively. The
horizontal line represents the imaginary axis in Figure $9$ and the line
$\Re(s)=-1$ in Figure $10$ respectively. The vertical one indicates the
real axis. The variation of the logarithmic scale is represented by gray
scale. If the absolute value is larger than $1.5$, the color of the
point is white, and the color of the point at which the absolute value
is the smallest (locally) is black. For the sake of a proper
understanding, we give also $3$-dimensional graphs for Figure $9$ and
$10$ in Figure $11$ and $12$ respectively. Notice that the height in
these $3$-dimensional figures is taken by an absolute value, not the
logarithm of it. We can see that the most left black holes of Figure $9$
corresponds to the trajectory of the non-trivial zero
$\r_1=\frac{1}{2}+\sqrt{-1}\times 14.13472\ldots$ and the most left
large hole soaked with black and the second black hole of Figure $10$
correspond to the trajectories of the trivial zero $s_1=-2$ and of
$\r_2=\frac{1}{2}+\sqrt{-1}\times 21.02203\ldots$ by looking at the
consequences of numerical calculations performing from $q=2^{-13}
\approx 1.2 \times 10^{-4}$ to $q=2^{-64}$. Also, the second black hole
in Figure $9$ seems to represent the point in the trajectory of $\r_5$ 
(see Figure $8$ and Remark\,\ref{rem:rho5}). It seems that there are $6$
and $7$ black holes which may indicate the zeros of $\z_q(s)$ and these
black holes lie on a line parallel to the imaginary axis. This comes
from the periodicity of the function $\frac{q^{s+r-1}}{1-q^{s+r-1}}$,
that is, $|\z_q(s+\d n)|$ is very small for $n=\pm 1,\pm 2,\ldots$
provided $\z_q(s)=0$. Actually, if $q$ is very small, then $\binom{s+\d
n+r-1}{r} \approx \binom{s+r-1}{r}$ and $(1-q)^{\d n} \approx
1$. Therefore, by Proposition\, \ref{prop:specialvalue}, we have
$\z_q(s+\d n)=(1-q)^{s+\d n}\sum^{\infty}_{r=0}\binom{s+\d n+r-1}{r}
\frac{q^{s+\d n+r-1}}{1-q^{s+\d n+r-1}} \approx \z_q(s)=0$. Thus, when
$q=2^{-64}$, i.e. $2\pi/{\log{q}}=-0.1416\ldots$, such a black hole
appears approximately every $0.14$ along the horizontal lines in Figure
$9$ and Figure $10$ respectively.

\begin{figure}[htbp]
\psfrag{Re}{$\text{Re}$}
\psfrag{Im}{$\text{Im}$}
\includegraphics[clip,width=162mm]{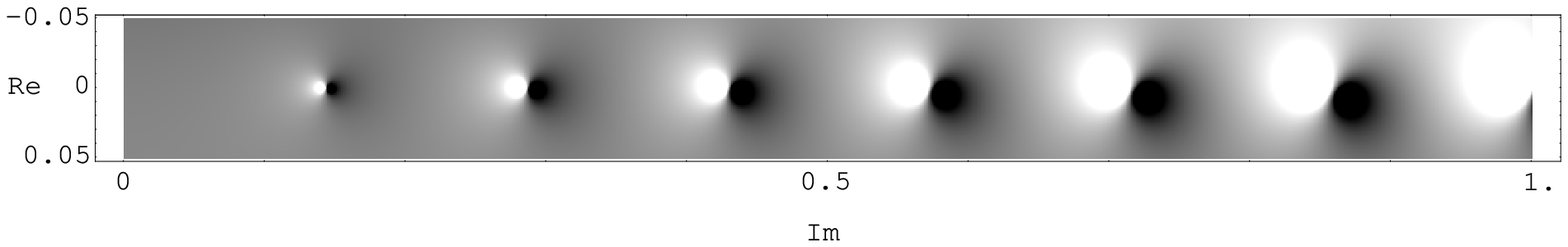}
\caption{$|\z_q(s)|\ \text{around $s=0$ at}\ q=2^{-64}\ \approx\ 5\times
 10^{-20}$} 
\end{figure}

\begin{figure}[htbp]
\psfrag{Re}{$\text{Re}$}
\psfrag{Im}{$\text{Im}$}
\includegraphics[clip,width=162mm]{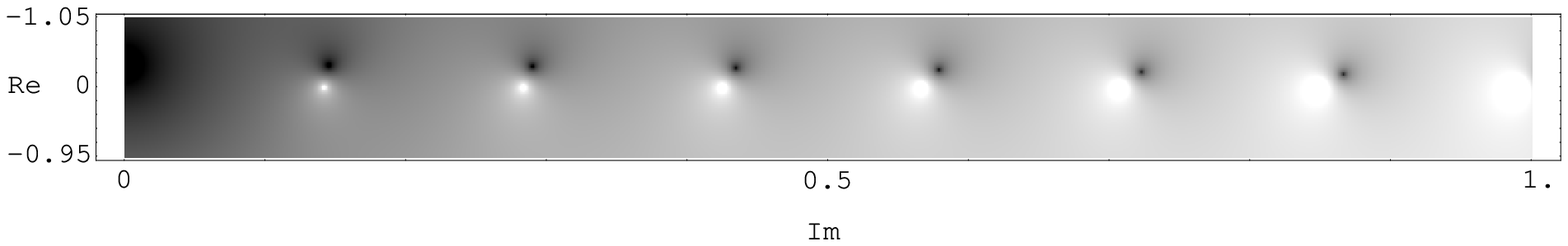}
\caption{$|\z_q(s)|\ \text{around $s=-1$ at}\ q=2^{-64}\ \approx\ 5\times
 10^{-20}$} 
\end{figure}

\begin{figure}[htbp]
\psfrag{x}{$\text{Re}$}
\psfrag{y}{$\text{Im}$}
\psfrag{Abs}{$|\z_q(s)|$}
\begin{center}
\includegraphics[clip]{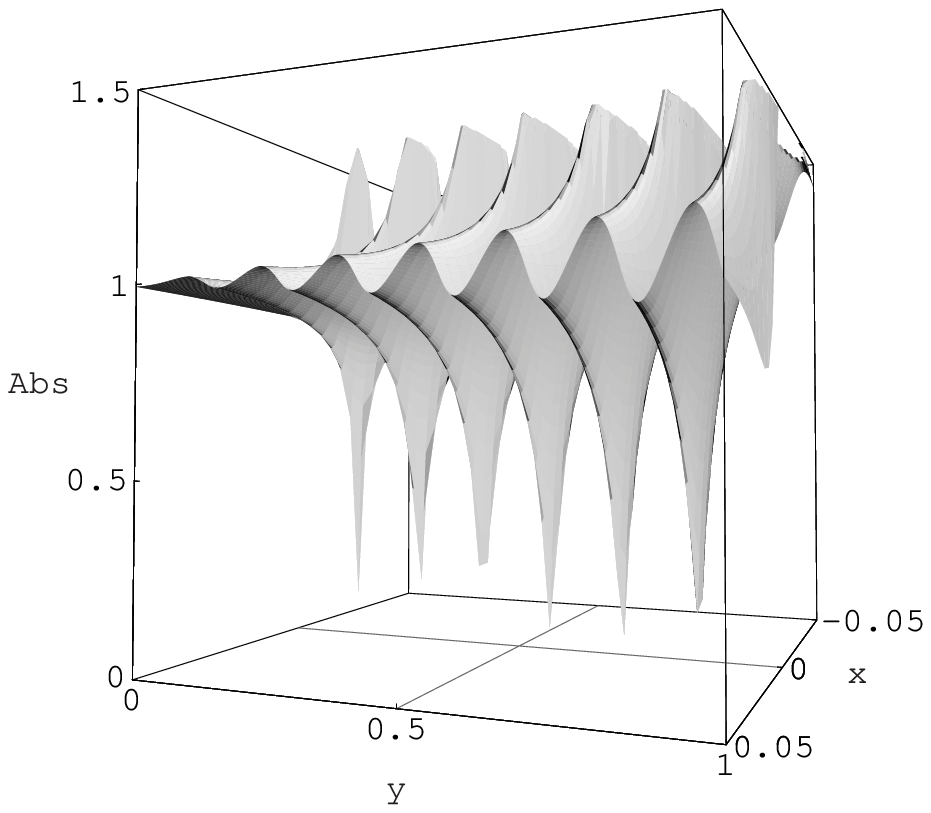}
\caption{$3$-dimensional graph for Figure\, $9$}
\end{center}
\end{figure}

\begin{figure}[htbp]
\psfrag{x}{$\text{Re}$}
\psfrag{y}{$\text{Im}$}
\psfrag{Abs}{$|\z_q(s)|$}
\begin{center}
\includegraphics[clip]{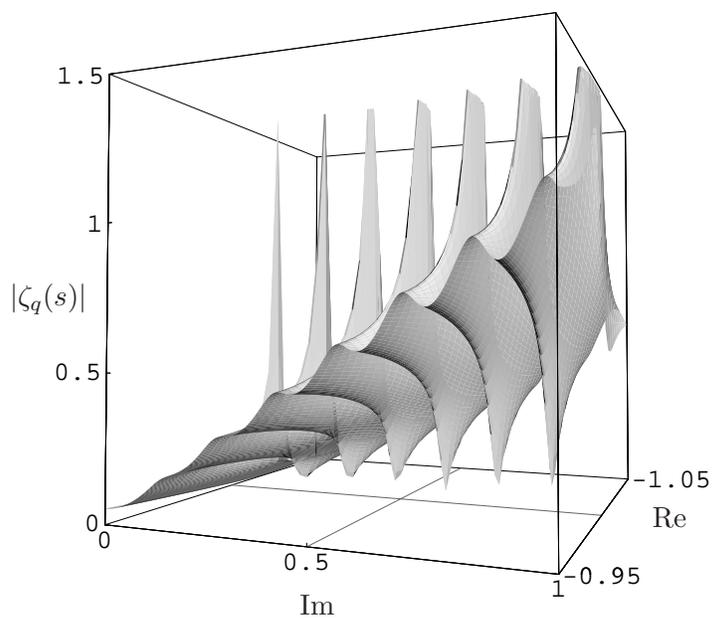}
\caption{$3$-dimensional graph for Figure\, $10$}
\end{center}
\end{figure}

 Proposition\,\ref{prop:tangent} indicates that the trajectory $z\nn(q)$
 can be tangent to the real axis when $q$ approaches to zero, but it is
 difficult to determine the approaching direction theoretically. From
 Figure $5$, $6$ and $7$, however, we can verify that the one of the
 points $\r^{(1)}_2(q)$, $\r^{(1)}_3(q)$, $\r^{(2)}_3(q)$ and
 $\r\nn_4(q)$ ($\n=1,2,3$) approaching a negative integer $-m$ seems to
 be tangent to the real axis from the left direction, while $\r\nn_1(q)$
 ($\n=1,2,3$), $\r^{(2)}_2(q)$, $\r^{(3)}_2(q)$, $\r^{(3)}_3(q)$ and
 $\r\nn_5(q)$ ($\n=1,2,3$) are not the cases. For instance, it seems
 that $\r^{(1)}_3(q)$ crosses the line $\Re(s)=-m$ first, then turns to
 the right and goes to $-m$ when $q$ approaches to $0$. This is
 consistent with the periodicity one can see in Figure $9$ and $10$
 mentioned above. 

 Lastly, we display the location of zeros $\r^{(1)}_j(q)$ of
 $\z_q(s)=\z^{(1)}_q(s)$ approximately which lie on the trajectory
 starting from the first seven non-trivial zeros of $\z(s)$ when $q=0.9$,
 $0.6$ and $0.3$ respectively in Figures $13$-$15$. In the figures, the
 number $j$ indicates $\r^{(1)}_j(q)$ for simplicity.

\begin{figure}[htbp]
 \centering
 \psfrag{Re}{$\Re$}
 \psfrag{Im}{$\Im$}
  \includegraphics[clip]{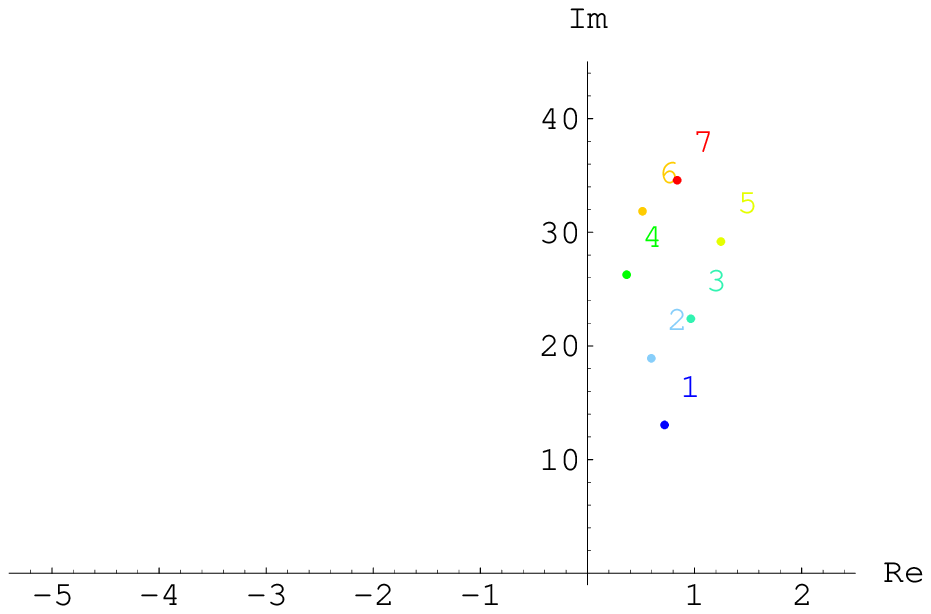}
 \caption{the zeros of $\z_q(s)$ correspond to the non-trivial zeros of  
 $\z(s)$ when $q=0.9$}  
\end{figure}

\begin{figure}[htbp]
 \centering
 \psfrag{Re}{$\Re$}
 \psfrag{Im}{$\Im$}
  \includegraphics[clip]{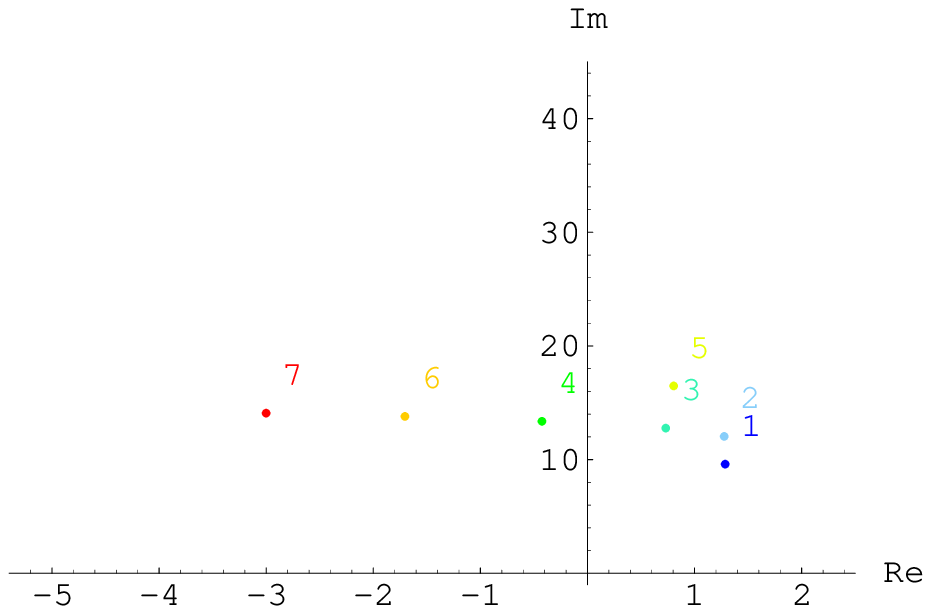}
 \caption{the zeros of $\z_q(s)$ correspond to the non-trivial zeros of
 $\z(s)$ when $q=0.6$}  
\end{figure}

\begin{figure}[htbp]
 \centering
 \psfrag{Re}{$\Re$}
 \psfrag{Im}{$\Im$}
  \includegraphics[clip]{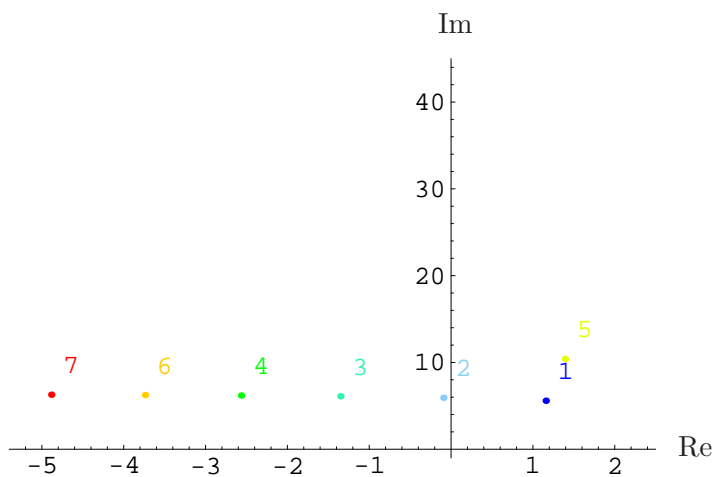}
 \caption{the zeros of $\z_q(s)$ correspond to the non-trivial zeros of
 $\z(s)$ when $q=0.3$}  
\end{figure}

\smallbreak
\smallbreak 
From the numerical calculation developed above, the following questions
naturally come up.

\smallbreak
\smallbreak


\noindent  
\textbf{Questions.}\textit{\ (i)\ Does the relation
$\Im(\r\nn_j(q))>\Im(\r\nn_l(q))\ (q\in(0,1))$ 
 hold for $j>l$? In particular, is $\r\nn_j(q)$ a simple zero of
 $\z\nn_q(s)$ whenever $\r\nn_j=\r\nn_j(1)$ is simple?\\   
 (ii)\ What is the shape of the function $\r\nn_j(q)$? For instance,
 where is the extremal point of the graph of $\Re(\r\nn_j(q))$ as a
 function of $\Im(\r\nn_j(q))$? etc. Moreover, what is an analogue of the
 Riemann hypothesis for $\z\nn_q(s)$?}\\      
 

\noindent
\textit{Acknowledgement.}
 We wish to thank Nobushige Kurokawa and Alberto Parmeggiani for their
 strong interest on the crystal limit behavior of the zeta functions and 
 related topics presented here.


\smallskip

\textsc{Kenichi Kawagoe}\\
Department of Computational Science, Faculty of Science, Kanazawa University, 
Kakuma Kanazawa 920-1192, Japan.\\
e-mail : \texttt{kawagoe@kenroku.kanazawa-u.ac.jp}\\

\textsc{Masato Wakayama}\\
Faculty of Mathematics, Kyushu University,
Hakozaki Fukuoka 812-8581, Japan.\\
e-mail : \texttt{wakayama@math.kyushu-u.ac.jp}\\

\textsc{Yoshinori Yamasaki}\\
Graduate School of Mathematics, Kyushu University,
Hakozaki Fukuoka 812-8581, Japan.\\
e-mail : \texttt{ma203032@math.kyushu-u.ac.jp}

\end{document}